\documentclass[11pt]{article}

\usepackage{amsthm,amsmath,amssymb,tikz,graphicx,array,arydshln,geometry}
\usepackage[colorlinks=true,
            linkcolor=blue,
            urlcolor=blue,
            citecolor=blue]{hyperref}

\newtheorem{thm}{Theorem}
\newtheorem{lem}[thm]{Lemma}

\newtheorem{prop}[thm]{Proposition}

\newtheorem{pr}[thm]{Question}

\theoremstyle{definition}
\newtheorem{defn}[thm]{Definition}

\theoremstyle{remark}

\newtheorem*{rem*}{Remark}

\tikzset{>=latex}
\usetikzlibrary{fit}
\tikzset{highlight/.style={thin,inner sep=0pt,rectangle,scale=1,rounded corners=1pt,draw}}
\newcommand{\tikzmark}[2]{%
\tikz[overlay,remember picture,baseline=(#1.base),inner sep=1.4pt] \node (#1) {#2};}
\newcommand{\HighlightEntry}[1][submatrix]{%
    \tikz[overlay,remember picture]{
    \node[highlight,fit=(#1)] {};}
}

\newcommand{\set}[1]{\big\{#1\big\}}
\newcommand{\floor}[1]{\left\lfloor #1\right\rfloor}
\newcommand{\ceil}[1]{\left\lceil #1\right\rceil}

\def\le{\leqslant}
\def\ge{\geqslant}

\def\ble{\prec_{B}}
\def\bleq{\preceq_{B}}
\def\sble{\prec_{\widehat{B}}}
\def\sbleq{\preceq_{\widehat{B}}}

\def\etal{~{\it et al.}}

\begin{document}

\title{On Maximum Chains in the Bruhat Order of $\mathcal{A}(n,2)$}

\author{M. Ghebleh\\[1ex]
{\em Department of Mathematics,
             Faculty of Science,
             Kuwait University,
             Kuwait}\\
\tt mohammad.ghebleh@ku.edu.kw}

\date{22 December 2013}

\maketitle

\begin{abstract}
Let $\mathcal{A}(R,S)$ denote the class of all matrices of zeros and ones with
row sum vector $R$ and column sum vector~$S$. We introduce the notion
of an inversion in a $(0,1)$--matrix.
This definition extends the standard notion of an inversion of a permutation,
in the sense that both notions agree on the class of permutation matrices.
We prove that the number of inversions in a $(0,1)$--matrix is monotonic
with respect to the secondary Bruhat order of the class $\mathcal{A}(R,S)$.
We apply this result in establishing
the maximum length of a chain in the Bruhat order of the class $\mathcal{A}(n,2)$
of $(0,1)$--matrices of order $n$ in which every row and every column has a sum of~$2$.
We give algorithmic constructions of chains of maximum
length in the Bruhat order of $\mathcal{A}(n,2)$.
\end{abstract}


\section{Introduction}

Let $R=(r_1,\ldots,r_m)$ and $S=(s_1,\ldots,s_n)$ be two vectors with
nonnegative integral entries.
The class of all $(0,1)$--matrices of size $m$ by $n$ with row sum vector $R$
and column sum vector $S$ is denoted by $\mathcal{A}(R,S)$.
If $m=n$ and $R=S=(k,k,\ldots,k)$, we simply write $\mathcal{A}(n,k)$ for $\mathcal{A}(R,S)$.
In particular, $\mathcal{A}(n,1)$ is the class of all permutation matrices of order $n$,
which can be identified by the symmetric group $\mathcal{S}_n$.
Combinatorial properties of the class $\mathcal{A}(R,S)$ are studied extensively
(see for example \cite{Brualdi1980,Brualdi2006,BrualdiBook2006,Ryser1963}
and the references there in).

Given the vectors $R,S$ and a matrix $A\in\mathcal{A}(R,S)$,
one may construct a new matrix $B\in\mathcal{A}(R,S)$ from $A$
by means of an {\em interchange}
\[
I_2=\left[\begin{array}{cc}1&0\\0&1\end{array}\right]
\leftrightarrow
\left[\begin{array}{cc}0&1\\1&0\end{array}\right]=L_2
\]
that replaces a $2\times2$ submatrix of $A$ equal to $I_2$ (if one exists)
by $L_2$, or vice versa.
In the class $\mathcal{A}({n,1})$ of permutation matrices, these interchanges correspond to
transpositions in permutations.
A result due to Ryser~\cite{Ryser1957,Ryser1963} states that
for any two matrices in the class $\mathcal{A}(R,S)$, one can be obtained from
the other by a sequence of $I_2\leftrightarrow L_2$ interchanges.

In~\cite{BrualdiHwang2004} Brualdi and Hwang define a {\em Bruhat order} on
the class $\mathcal{A}(R,S)$ generalizing the classical Bruhat order on the symmetric
group $\mathcal{S}_n$ (note that the class $\mathcal{A}(n,1)$ consists of permutation matrices
of order $n$, hence it can be identified by $\mathcal{S}_n$).
Given a $(0,1)$--matrix $A$ of size $m$ by $n$, let $\Sigma_A$ be the $m$ by $n$
matrix whose $(k,\ell)$--entry is
\[
\sigma_{k\ell}(A)=\sum_{i=1}^{k}\sum_{j=1}^{\ell}a_{ij}.
\]
If $A$ and $C$ are $(0,1)$--matrices in a class $\mathcal{A}(R,S)$, then 
$A$ precedes $C$ in the Bruhat order, written as $A\bleq C$ for short,
if $\Sigma_A\ge\Sigma_C$ in the entrywise order. Namely,
\[
A\bleq C\text{ \ if and only if \ }\sigma_{ij}(A)\ge\sigma_{ij}(C)
\text{ \ for all \ }1\le i\le m\text{ \ and \ }1\le j\le n.
\]
It is easily observed that if $C$ is obtained from $A$ by a sequence of
$I_2\to L_2$ interchanges, then $A\bleq C$. It is shown in~\cite{BrualdiDeaett2007}
that the converse does not hold in general. This observation defines a
{\em secondary Bruhat order} on the classes $\mathcal{A}(R,S)$ of $(0,1)$--matrices:
$A\sbleq C$ if and only if $C$ is obtained from $A$ by a sequence of
$I_2\to L_2$ interchanges.
It is shown in~\cite{BrualdiDeaett2007} that the Bruhat order and the
secondary Bruhat order are the same on the classes $\mathcal{A}(n,2)$,
but they are different on $\mathcal{A}(6,3)$.

Answering a question asked in~\cite{BrualdiDeaett2007},
Conflitti\etal~\cite{ConflittiDaFonsecaMamede2012}
show that for all $k\ge1$, the maximum length of a chain in the
Bruhat order of the class $\mathcal{A}({2k,k})$ is~$k^4$.
In this work, we establish the maximum length of a chain in the Bruhat order of
the classes $\mathcal{A}(n,2)$. We define the notion of an {inversion} in a
$(0,1)$--matrix and show that the number of inversions in a $(0,1)$--matrix
is monotonic with respect
to the secondary Bruhat order. This result, together with a classification of
minimal elements of the Bruhat order of $\mathcal{A}(n,2)$ proved in~\cite{BrualdiDeaett2007}, gives an upper bound on the
length of a chain in the (secondary) Bruhat order of the class $\mathcal{A}(n,2)$.
We give algorithmic constructions of chains that achieve this upper bound.

\section{Inversions in $(0,1)$--matrices}

The symmetric group $\mathcal{S}_n$ is naturally identified with the class $\mathcal{A}(n,1)$
of permutation matrices of order~$n$. In this sense,
an inversion in a permutation corresponds to a pair of ones in
the corresponding permutation matrix $P$, one of which is located
to the top-right of the other. More precisely, an inversion in $P$ consists
of two ones in the entries $(i,j)$ and $(i',j')$, such that $i<i'$ and $j>j'$.
We adopt the same definition for an inversion in any $(0,1)$--matrix.

\begin{defn}
Let $A=[a_{ij}]$ be a $(0,1)$--matrix. An {\em inversion} in $A$ consists of
any two entries $a_{ij}=a_{k\ell}=1$ such that $(i-k)(j-\ell)<0$.
We denote the total number of inversions in $A$ by $\nu(A)$.
\end{defn}

This definition is illustrated in the following matrix $A$, where each
of its $\nu(A)=9$ inversions is represented by a line segment between
the two entries it involves.
\[
\def\arraystretch{1.2}
\def\tab{1.4em}
A=\left[\begin{array}{c@{\hspace{\tab}}c@{\hspace{\tab}}c@{\hspace{\tab}}c@{\hspace{\tab}}c}
1 & \tikzmark{a12}{1} & \tikzmark{a13}{1} & 0 & \tikzmark{a15}{1}\\
\tikzmark{a21}{1} & 0 & 0 & 0 & 0\\
0 & \tikzmark{a32}{1} & 0 & 0 & \tikzmark{a35}{1}\\
0 & 0 & \tikzmark{a43}{1} & \tikzmark{a44}{1} & 0
\end{array}\right]
\]
\HighlightEntry[a12]
\HighlightEntry[a13]
\HighlightEntry[a21]
\HighlightEntry[a15]
\HighlightEntry[a32]
\HighlightEntry[a35]
\HighlightEntry[a43]
\HighlightEntry[a44]
\tikz[overlay,remember picture]{
\draw (a12)--(a21);
\draw (a13)--(a21);
\draw (a15)--(a21);
\draw (a13)--(a32);
\draw (a15)--(a32);
\draw (a15)--(a43);
\draw (a15)--(a44);
\draw (a35)--(a43);
\draw (a35)--(a44);
}

\begin{lem}
Let $A,C\in\mathcal{A}(R,S)$. If $A\sble C$, then $\nu(A)<\nu(C)$.
\label{lem:nu}
\end{lem}

\begin{proof}
By definition, there is a sequence of $I_2\to L_2$ interchanges
which transforms $A$ to~$C$. Thus it suffices to prove that the number
of inversions increases under each interchange.
Let $D$ be a $(0,1)$--matrix such that $D[\{i,i'\},\{j,j'\}]=I_2$.
That is, the submatrix of $D$ at the intersection of rows
$i$ and $i'$ and columns $j$ and $j'$ equals~$I_2$.
Let $a,b,c,d,e$ denote the number of ones in different blocks of $D$
as shown in the following.
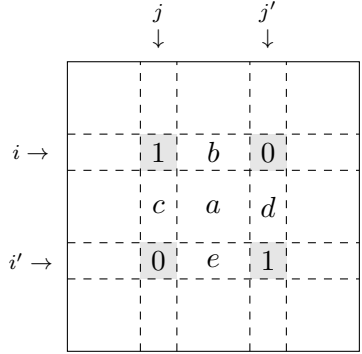
\begin{center}
\vskip -1ex
\begin{tikzpicture}[scale=0.48]
\draw (0,0) rectangle (8,8);
\foreach \x in {2,5} \foreach \y in {2,5}
	\fill [black!10!white] (\x,\y) rectangle (1+\x,1+\y);
\foreach \x in {2,3,5,6} {
	\draw[dashed] (\x,0)--(\x,8);
	\draw[dashed] (0,\x)--(8,\x);
}
\node at (2.5,2.5) {$0$};
\node at (5.5,5.5) {$0$};
\node at (2.5,5.5) {$1$};
\node at (5.5,2.5) {$1$};
\node at (4,5.5) {$b$};
\node at (2.5,4) {$c$};
\node at (4,4) {$a$};
\node at (5.5,4) {$d$};
\node at (4,2.5) {$e$};
\node at (-1,5.5) {\scriptsize $i\to$};
\node at (9.5,5.5) {\ };
\node at (-1,2.5) {\scriptsize $i'\to$};
\node at (2.5,9) {\scriptsize\begin{tabular}{c}$j$\\$\downarrow$\end{tabular}};
\node at (5.5,9) {\scriptsize\begin{tabular}{c}$j'$\\$\downarrow$\end{tabular}};
\end{tikzpicture}
\end{center}
Let $D'$ be obtained from $D$ by interchanging this $I_2$ to $L_2$.
By counting the number of inversions in $D$ that involve either of
the ones at the $(i,j)$ and $(i',j')$ entries, and similarly the number of
inversions in $D'$ that involve either of the ones at the $(i,j')$ and
$(i',j)$ entries, it is easy to see that
$$\nu(D')=\nu(D)+1+2a+b+c+d+e\ge\nu(D)+1.\qedhere$$
\end{proof}

Brualdi and Deaett prove in~\cite{BrualdiDeaett2007} that for $\mathcal{A}(n,2)$, the Bruhat order and the secondary Bruhat order are the same. The next lemma is an immediate corollary of this result and Lemma~\ref{lem:nu}.

\begin{lem}
Let $A,C\in\mathcal{A}(n,2)$. If $A\ble C$, then $\nu(A)<\nu(C)$.
\label{lem:nuAn2}
\end{lem}

Note that $\nu(A)<\nu(C)$ does not necessarily imply $A\ble C$. For example, if
\[
A=\left[\begin{array}{cccc}1&0&0&1\\1&1&0&0\\0&1&1&0\\0&0&1&1\end{array}\right]
\text{ \ and \ }
C=\left[\begin{array}{cccc}0&1&1&0\\1&1&0&0\\1&0&0&1\\0&0&1&1\end{array}\right],
\]
then 
$\nu(A)=5$ and $\nu(C)=7$, while
\[\sigma_{11}(A)=1>0=\sigma_{11}(C)
\text{ \ and \ }
\sigma_{13}(A)=1<2=\sigma_{13}(C),\]
showing that $A$ and $C$ are incomparable in the Bruhat order of $\mathcal{A}(4,2)$.

In~\cite{BrualdiDeaett2007} two matrices in $\mathcal{A}(6,3)$ are given that are
comparable in the Bruhat order and incomparable in the secondary Bruhat order, thus proving that these two orders are not the same on $\mathcal{A}(6,3)$. The result in Lemma~\ref{lem:nuAn2} on the other hand, seems to hold for all classes $\mathcal{A}(R,S)$.
Our computational experiments, while far from being thorough,
point to an affirmative answer to the following question.

\begin{pr}
Let $A,C\in\mathcal{A}(R,S)$. Does $A\ble C$ imply $\nu(A)<\nu(C)$?
\label{q:ARSmonotone}
\end{pr}

\section{Maximum chains in the Bruhat order of $\mathcal{A}(n,2)$}

Any maximal chain in a partially ordered set begins with a minimal element and ends with a maximal element, since otherwise, it could be extended to a larger chain.
Minimal matrices in the Bruhat order of the class $\mathcal{A}(n,2)$ are characterized in~\cite{BrualdiHwang2004}. It is shown there that a matrix in $\mathcal{A}(n,2)$ is minimal in the Bruhat order if and only if it is the direct sum of matrices each of which equals
\[
J_2=\left[\begin{array}{cc}1&1\\1&1\end{array}\right]
\text{ \ or \ }
F_3=\left[\begin{array}{ccc}1&1&0\\1&0&1\\0&1&1\end{array}\right].
\]
On the other hand, if $A,C\in\mathcal{A}(n,2)$ such that $A\bleq C$, and $A'$ and $C'$ are obtained from $A$ and $C$ respectively by reversing the order of their columns (flipping the matrix in the left/right direction), then $C'\bleq A'$.
Therefore, maximal matrices in the Bruhat order of $\mathcal{A}(n,2)$ are obtained by reversing the order of columns in minimal matrices. For every $n\ge4$, we consider one special minimal matrix $P_n$ in the Bruhat order of
$\mathcal{A}(n,2)$, and a maximal matrix $Q_n$ that is obtained from $P_n$ by reversing
the order of its columns. If $n$ is even, $P_n$ is the direct sum of $n/2$ copies of~$J_2$:
\[
\def\arraystretch{1.2}
P_n=\left[\begin{array}{cccc}
J_2&0&\cdots&0\\
0&J_2&\cdots&0\\
\vdots&\vdots & \ddots&\vdots\\
0&0&\cdots&J_2
\end{array}\right]
\text{ \ \ and \ \ }
Q_n=\left[\begin{array}{cccc}
0&\cdots&0&J_2\\
0&\cdots&J_2&0\\
\vdots & \reflectbox{$\ddots$}&\vdots&\vdots\\
J_2&\cdots&0&0
\end{array}\right],
\]
and if $n$ is odd, $P_n$ is the direct sum of $(n-3)/2$ copies of $J_2$ and one copy of $F_3$:
\[
\def\arraystretch{1.2}
P_n=\left[\begin{array}{ccccc}
J_2&0&\cdots&0&0\\
0&J_2&\cdots&0&0\\
\vdots&\vdots & \ddots&\vdots&\vdots\\
0&0&\cdots&J_2&0\\
0&0&\cdots&0&F_3
\end{array}\right]
\text{ \ \ and \ \ }
Q_n=\left[\begin{array}{ccccc}
0&0&\cdots&0&J_2\\
0&0&\cdots&J_2&0\\
\vdots &\vdots & \reflectbox{$\ddots$}&\vdots&\vdots\\
0&J_2&\cdots&0&0\\
F_3'&0&\cdots&0&0
\end{array}\right].
\]

In this section, we first give constructions of chains from $P_n$ to $Q_n$
in the Bruhat order of $\mathcal{A}(n,2)$.
We then prove in Theorem~\ref{thm:maxlen}
that these chains are indeed the longest possible chains in the Bruhat order of $\mathcal{A}(n,2)$.
The next three lemmas, the first of which is a special case of
the main result of~\cite{ConflittiDaFonsecaMamede2012},
provide the necessary ingredients in our constructions.

\begin{lem}{\rm\cite{ConflittiDaFonsecaMamede2012}}
There is a chain of length $16$ from $P_4$ to $Q_4$ in the Bruhat order
of $\mathcal{A}(4,2)$.
\label{lem:chain16}
\end{lem}

\begin{lem}
There is a chain of length $29$ from $P_5$ to $Q_5$ in the Bruhat order
of $\mathcal{A}(5,2)$.
\label{lem:chain29}
\end{lem}

\begin{proof}
A chain of length $6$ from $P_5$ to the matrix
\[
Z=\left[\begin{array}{ccccc}
1 & 1 & 0 & 0 & 0 \\
1 & 0 & 0 & 1 & 0 \\
0 & 1 & 0 & 0 & 1 \\
0 & 0 & 1 & 0 & 1 \\
0 & 0 & 1 & 1 & 0 \\
\end{array}\right].
\]
is presented in Figure~\ref{fig:chain6}
and a chain of length $23$ from $Z$ to $Q_5$ is presented in Figure~\ref{fig:chain23}.
The desired chain is the concatenation of these two chains.
\end{proof}

\begin{figure}[ht]
\begin{center}
\def\arraystretch{.84}
$\begin{array}{@{}c@{\ }c@{\ }c@{\ }c@{\ }c@{\ }c@{\ }c@{\ }c@{}}
&
\left[\begin{array}{@{\,}c@{\ \,}c@{\ \,}c@{\ \,}c@{\ \,}c@{\,}}1&1&\cdot&\cdot&\cdot\\1&1&\cdot&\cdot&\cdot\\\cdot&\cdot&1&1&\cdot\\\cdot&\cdot&1&\cdot&1\\\cdot&\cdot&\cdot&1&1\\\end{array}\right]
& \to &
\left[\begin{array}{@{\,}c@{\ \,}c@{\ \,}c@{\ \,}c@{\ \,}c@{\,}}1&1&\cdot&\cdot&\cdot\\1&1&\cdot&\cdot&\cdot\\\cdot&\cdot&1&1&\cdot\\\cdot&\cdot&\cdot&1&1\\\cdot&\cdot&1&\cdot&1\\\end{array}\right]
& \to &
\left[\begin{array}{@{\,}c@{\ \,}c@{\ \,}c@{\ \,}c@{\ \,}c@{\,}}1&1&\cdot&\cdot&\cdot\\1&\cdot&1&\cdot&\cdot\\\cdot&1&\cdot&1&\cdot\\\cdot&\cdot&\cdot&1&1\\\cdot&\cdot&1&\cdot&1\\\end{array}\right]
& \to &
\left[\begin{array}{@{\,}c@{\ \,}c@{\ \,}c@{\ \,}c@{\ \,}c@{\,}}1&1&\cdot&\cdot&\cdot\\1&\cdot&\cdot&1&\cdot\\\cdot&1&1&\cdot&\cdot\\\cdot&\cdot&\cdot&1&1\\\cdot&\cdot&1&\cdot&1\\\end{array}\right]
\\ \\ \to &
\left[\begin{array}{@{\,}c@{\ \,}c@{\ \,}c@{\ \,}c@{\ \,}c@{\,}}1&1&\cdot&\cdot&\cdot\\1&\cdot&\cdot&1&\cdot\\\cdot&1&\cdot&1&\cdot\\\cdot&\cdot&1&\cdot&1\\\cdot&\cdot&1&\cdot&1\\\end{array}\right]
& \to &
\left[\begin{array}{@{\,}c@{\ \,}c@{\ \,}c@{\ \,}c@{\ \,}c@{\,}}1&1&\cdot&\cdot&\cdot\\1&\cdot&\cdot&1&\cdot\\\cdot&1&\cdot&\cdot&1\\\cdot&\cdot&1&1&\cdot\\\cdot&\cdot&1&\cdot&1\\\end{array}\right]
& \to &
\left[\begin{array}{@{\,}c@{\ \,}c@{\ \,}c@{\ \,}c@{\ \,}c@{\,}}1&1&\cdot&\cdot&\cdot\\1&\cdot&\cdot&1&\cdot\\\cdot&1&\cdot&\cdot&1\\\cdot&\cdot&1&\cdot&1\\\cdot&\cdot&1&1&\cdot\\\end{array}\right]
\end{array}$
\vspace{-3.2ex}
\end{center}
\caption{A chain of length $6$ from $P_5$ to $Z$ in the Bruhat order of $\mathcal{A}(5,2)$. Dots represent zeros.\label{fig:chain6}}
\end{figure}

\begin{figure}[ht]
\begin{center}
\def\arraystretch{.84}
$\begin{array}{@{}c@{\ }c@{\ }c@{\ }c@{\ }c@{\ }c@{\ }c@{\ }c@{\ }c@{\ }c@{}}
&\left[\begin{array}{@{\,}c@{\ \,}c@{\ \,}c@{\ \,}c@{\ \,}c@{\,}}1&1&\cdot&\cdot&\cdot\\1&\cdot&\cdot&1&\cdot\\\cdot&1&\cdot&\cdot&1\\\cdot&\cdot&1&\cdot&1\\\cdot&\cdot&1&1&\cdot\\\end{array}\right]
& \to & 
\left[\begin{array}{@{\,}c@{\ \,}c@{\ \,}c@{\ \,}c@{\ \,}c@{\,}}1&1&\cdot&\cdot&\cdot\\1&\cdot&\cdot&\cdot&1\\\cdot&1&\cdot&1&\cdot\\\cdot&\cdot&1&\cdot&1\\\cdot&\cdot&1&1&\cdot\\\end{array}\right]
& \to & 
\left[\begin{array}{@{\,}c@{\ \,}c@{\ \,}c@{\ \,}c@{\ \,}c@{\,}}1&1&\cdot&\cdot&\cdot\\1&\cdot&\cdot&\cdot&1\\\cdot&\cdot&1&1&\cdot\\\cdot&1&\cdot&\cdot&1\\\cdot&\cdot&1&1&\cdot\\\end{array}\right]
& \to & 
\left[\begin{array}{@{\,}c@{\ \,}c@{\ \,}c@{\ \,}c@{\ \,}c@{\,}}1&1&\cdot&\cdot&\cdot\\1&\cdot&\cdot&\cdot&1\\\cdot&\cdot&1&1&\cdot\\\cdot&\cdot&1&\cdot&1\\\cdot&1&\cdot&1&\cdot\\\end{array}\right]
\\ \\ \to & 
\left[\begin{array}{@{\,}c@{\ \,}c@{\ \,}c@{\ \,}c@{\ \,}c@{\,}}1&1&\cdot&\cdot&\cdot\\1&\cdot&\cdot&\cdot&1\\\cdot&\cdot&1&1&\cdot\\\cdot&\cdot&\cdot&1&1\\\cdot&1&1&\cdot&\cdot\\\end{array}\right]
& \to & 
\left[\begin{array}{@{\,}c@{\ \,}c@{\ \,}c@{\ \,}c@{\ \,}c@{\,}}1&\cdot&1&\cdot&\cdot\\1&\cdot&\cdot&\cdot&1\\\cdot&1&\cdot&1&\cdot\\\cdot&\cdot&\cdot&1&1\\\cdot&1&1&\cdot&\cdot\\\end{array}\right]
& \to &
\left[\begin{array}{@{\,}c@{\ \,}c@{\ \,}c@{\ \,}c@{\ \,}c@{\,}}1&\cdot&\cdot&1&\cdot\\1&\cdot&\cdot&\cdot&1\\\cdot&1&1&\cdot&\cdot\\\cdot&\cdot&\cdot&1&1\\\cdot&1&1&\cdot&\cdot\\\end{array}\right]
& \to & 
\left[\begin{array}{@{\,}c@{\ \,}c@{\ \,}c@{\ \,}c@{\ \,}c@{\,}}1&\cdot&\cdot&1&\cdot\\1&\cdot&\cdot&\cdot&1\\\cdot&1&\cdot&1&\cdot\\\cdot&\cdot&1&\cdot&1\\\cdot&1&1&\cdot&\cdot\\\end{array}\right]
\\ \\ \to & 
\left[\begin{array}{@{\,}c@{\ \,}c@{\ \,}c@{\ \,}c@{\ \,}c@{\,}}1&\cdot&\cdot&1&\cdot\\1&\cdot&\cdot&\cdot&1\\\cdot&1&\cdot&\cdot&1\\\cdot&\cdot&1&1&\cdot\\\cdot&1&1&\cdot&\cdot\\\end{array}\right]
& \to & 
\left[\begin{array}{@{\,}c@{\ \,}c@{\ \,}c@{\ \,}c@{\ \,}c@{\,}}1&\cdot&\cdot&\cdot&1\\1&\cdot&\cdot&1&\cdot\\\cdot&1&\cdot&\cdot&1\\\cdot&\cdot&1&1&\cdot\\\cdot&1&1&\cdot&\cdot\\\end{array}\right]
& \to & 
\left[\begin{array}{@{\,}c@{\ \,}c@{\ \,}c@{\ \,}c@{\ \,}c@{\,}}1&\cdot&\cdot&\cdot&1\\1&\cdot&\cdot&1&\cdot\\\cdot&\cdot&1&\cdot&1\\\cdot&1&\cdot&1&\cdot\\\cdot&1&1&\cdot&\cdot\\\end{array}\right]
& \to &
\left[\begin{array}{@{\,}c@{\ \,}c@{\ \,}c@{\ \,}c@{\ \,}c@{\,}}1&\cdot&\cdot&\cdot&1\\\cdot&1&\cdot&1&\cdot\\\cdot&\cdot&1&\cdot&1\\1&\cdot&\cdot&1&\cdot\\\cdot&1&1&\cdot&\cdot\\\end{array}\right]
\\ \\ \to & 
\left[\begin{array}{@{\,}c@{\ \,}c@{\ \,}c@{\ \,}c@{\ \,}c@{\,}}1&\cdot&\cdot&\cdot&1\\\cdot&1&\cdot&1&\cdot\\\cdot&\cdot&1&\cdot&1\\\cdot&1&\cdot&1&\cdot\\1&\cdot&1&\cdot&\cdot\\\end{array}\right]
& \to & 
\left[\begin{array}{@{\,}c@{\ \,}c@{\ \,}c@{\ \,}c@{\ \,}c@{\,}}1&\cdot&\cdot&\cdot&1\\\cdot&1&\cdot&1&\cdot\\\cdot&\cdot&1&\cdot&1\\\cdot&\cdot&1&1&\cdot\\1&1&\cdot&\cdot&\cdot\\\end{array}\right]
& \to & 
\left[\begin{array}{@{\,}c@{\ \,}c@{\ \,}c@{\ \,}c@{\ \,}c@{\,}}1&\cdot&\cdot&\cdot&1\\\cdot&\cdot&1&1&\cdot\\\cdot&1&\cdot&\cdot&1\\\cdot&\cdot&1&1&\cdot\\1&1&\cdot&\cdot&\cdot\\\end{array}\right]
& \to & 
\left[\begin{array}{@{\,}c@{\ \,}c@{\ \,}c@{\ \,}c@{\ \,}c@{\,}}1&\cdot&\cdot&\cdot&1\\\cdot&\cdot&1&1&\cdot\\\cdot&\cdot&1&\cdot&1\\\cdot&1&\cdot&1&\cdot\\1&1&\cdot&\cdot&\cdot\\\end{array}\right]
\\ \\ \to &
\left[\begin{array}{@{\,}c@{\ \,}c@{\ \,}c@{\ \,}c@{\ \,}c@{\,}}1&\cdot&\cdot&\cdot&1\\\cdot&\cdot&1&1&\cdot\\\cdot&\cdot&\cdot&1&1\\\cdot&1&1&\cdot&\cdot\\1&1&\cdot&\cdot&\cdot\\\end{array}\right]
& \to & 
\left[\begin{array}{@{\,}c@{\ \,}c@{\ \,}c@{\ \,}c@{\ \,}c@{\,}}\cdot&\cdot&1&\cdot&1\\1&\cdot&\cdot&1&\cdot\\\cdot&\cdot&\cdot&1&1\\\cdot&1&1&\cdot&\cdot\\1&1&\cdot&\cdot&\cdot\\\end{array}\right]
& \to & 
\left[\begin{array}{@{\,}c@{\ \,}c@{\ \,}c@{\ \,}c@{\ \,}c@{\,}}\cdot&\cdot&\cdot&1&1\\1&\cdot&1&\cdot&\cdot\\\cdot&\cdot&\cdot&1&1\\\cdot&1&1&\cdot&\cdot\\1&1&\cdot&\cdot&\cdot\\\end{array}\right]
& \to & 
\left[\begin{array}{@{\,}c@{\ \,}c@{\ \,}c@{\ \,}c@{\ \,}c@{\,}}\cdot&\cdot&\cdot&1&1\\1&\cdot&\cdot&1&\cdot\\\cdot&\cdot&1&\cdot&1\\\cdot&1&1&\cdot&\cdot\\1&1&\cdot&\cdot&\cdot\\\end{array}\right]
\\ \\ \to & 
\left[\begin{array}{@{\,}c@{\ \,}c@{\ \,}c@{\ \,}c@{\ \,}c@{\,}}\cdot&\cdot&\cdot&1&1\\1&\cdot&\cdot&\cdot&1\\\cdot&\cdot&1&1&\cdot\\\cdot&1&1&\cdot&\cdot\\1&1&\cdot&\cdot&\cdot\\\end{array}\right]
& \to &
\left[\begin{array}{@{\,}c@{\ \,}c@{\ \,}c@{\ \,}c@{\ \,}c@{\,}}\cdot&\cdot&\cdot&1&1\\\cdot&\cdot&1&\cdot&1\\1&\cdot&\cdot&1&\cdot\\\cdot&1&1&\cdot&\cdot\\1&1&\cdot&\cdot&\cdot\\\end{array}\right]
& \to & 
\left[\begin{array}{@{\,}c@{\ \,}c@{\ \,}c@{\ \,}c@{\ \,}c@{\,}}\cdot&\cdot&\cdot&1&1\\\cdot&\cdot&\cdot&1&1\\1&\cdot&1&\cdot&\cdot\\\cdot&1&1&\cdot&\cdot\\1&1&\cdot&\cdot&\cdot\\\end{array}\right]
& \to & 
\left[\begin{array}{@{\,}c@{\ \,}c@{\ \,}c@{\ \,}c@{\ \,}c@{\,}}\cdot&\cdot&\cdot&1&1\\\cdot&\cdot&\cdot&1&1\\\cdot&1&1&\cdot&\cdot\\1&\cdot&1&\cdot&\cdot\\1&1&\cdot&\cdot&\cdot\\\end{array}\right]
\end{array}$
\vspace{-3.2ex}
\end{center}
\caption{A chain of length $23$ from $Z$ to $Q_5$ in the Bruhat order of $\mathcal{A}(5,2)$. Dots represent zeros.\label{fig:chain23}}
\end{figure}

\begin{lem}
Let $Y$ denote the direct sum of the matrices $J_2$ and $F_3'$.
There is a chain of length $24$ from $Y$ to $Q_5$ in the Bruhat order of
$\mathcal{A}({5,2})$.
\label{lem:chain24}
\end{lem}

\begin{proof}
Let $Z$ be as defined in the proof of Lemma~\ref{lem:chain29}. By computing
the matrices $\Sigma_Y$ and $\Sigma_Z$, it is easily observed that $Y\bleq Z$.
The desired chain is obtained by extending the chain of Figure~\ref{fig:chain23} at the beginning with~$Y$.
\end{proof}

Let $A_0,A_1,\ldots,A_k$ be $(0,1)$--matrices of the same size $m\times n$,
and let $r$ and $s$ be positive integers.
Let $1\le i_1<i_2<\cdots<i_r\le m$ and $1\le j_1<j_2<\cdots<j_s\le n$ be constant
indices.
Let $L_i=A[\{i_1,i_2,\ldots,i_r\},\{j_1,j_2,\ldots,j_s\}]$ be the $r\times s$
submatrix of $A_i$ consisting of all entries in
positions $(i_k,j_\ell)$ where $1\le k\le r$ and $1\le \ell\le s$.
Suppose that
$L_0\bleq L_1\bleq\cdots\bleq L_k$,
and that all matrices $A_i$ agree in entries that lie outside
the submatrices $L_i$. Then it is clear that
$A_0\bleq A_1\bleq\cdots\bleq A_k$.
This observation allows constructions of chains in the Bruhat order of $\mathcal{A}(n,2)$
using the chains presented in Lemmas~\ref{lem:chain16}, \ref{lem:chain29}, and~\ref{lem:chain24}.

\begin{prop}
If $n\ge4$ is even, then there is a chain of length $2n(n-2)$ from $P_n$ to $Q_n$ in the Bruhat order of $\mathcal{A}(n,2)$.
\label{prop:chainEven}
\end{prop}

\begin{proof}
All chains in this proof are with respect to the Bruhat order of $(0,1)$--matrices.
We give a recursive construction, thus proving the existence of the desired chain
by induction on~$n$.
For $n=4$ the chain of Lemma~\ref{lem:chain16} has the desired length.
Let $n\ge6$. 
Note that $P_n=P_{n-2}\oplus J_2$, where $\oplus$ denotes the direct sum of matrices.
By the induction hypothesis, there is a chain of length $2(n-2)(n-4)$
from $P_{n-2}$ to $Q_{n-2}$. Taking the direct sum of the matrices in such
chain with $J_2$, we obtain a chain of the same length from $P_n$ to
$A_1=Q_{n-2}\oplus J_2$.
We extend this chain to one from $P_n$ to $Q_n$ as follows.
Let $E_1$ be the submatrix of $A_1$ induced by rows $1,2,n-1,n$ and columns
$n-3,n-2,n-1,n$. Then $E_1=P_4$. We extend the current chain by 
keeping all entries outside $E_1$ constant, and applying the chain of
Lemma~\ref{lem:chain16} in the positions corresponding to $E_1$.
This extends the current chain by $16$. Let $A_2$ denote the end of this
chain.
We proceed by applying the same procedure to the submatrix $E_2=P_4$
of $A_2$ induced by rows $3,4,n-1,n$ and columns $n-5,n-4,n-3,n-2$.
We repeat this procedure for a total number of $n/2-1$ rounds,
after which the resulting chain ends at $A_{n/2}=Q_n$. The length of
this chain is
\[
2(n-2)(n-4)+16(n/2-1)=(n-2)(2n-8+8)=2n(n-2).
\]
The algorithm presented in this proof is illustrated in Figure~\ref{fig:proofEven} for $n=8$.
\end{proof}

\begin{figure}[ht]
\begin{center}
\scalebox{.96}{$\begin{array}{@{}c@{}c@{}c@{}c@{}c@{}c@{}c@{}}
\left[\begin{array}{@{}cccc@{}}
J&\cdot&\cdot&\cdot\\
\cdot&J&\cdot&\cdot\\
\cdot&\cdot&J&\cdot\\
\cdot&\cdot&\cdot&J\\
\end{array}\right]
&
\begin{array}{c}\scalebox{.8}{48}\\[-1ex]\longrightarrow\end{array}
&
\left[\begin{array}{@{}cccc@{}}
\cdot&\cdot&J&\cdot\\
\cdot&J&\cdot&\cdot\\
J&\cdot&\cdot&\cdot\\
\cdot&\cdot&\cdot&J\\
\end{array}\right]
&\multicolumn{4}{l}{\small\text{(Using the induction hypothesis for }n=6\text{)}}
\\ \\[-1em]
&
\begin{array}{c}\scalebox{.8}{16}\\[-1ex]\longrightarrow\end{array}
&
\left[\begin{array}{@{}cccc@{}}
\cdot&\cdot&\cdot&J\\
\cdot&J&\cdot&\cdot\\
J&\cdot&\cdot&\cdot\\
\cdot&\cdot&J&\cdot\\
\end{array}\right]
&
\begin{array}{c}\scalebox{.8}{16}\\[-1ex]\longrightarrow\end{array}
&
\left[\begin{array}{@{}cccc@{}}
\cdot&\cdot&\cdot&J\\
\cdot&\cdot&J&\cdot\\
J&\cdot&\cdot&\cdot\\
\cdot&J&\cdot&\cdot\\
\end{array}\right]
&
\begin{array}{c}\scalebox{.8}{16}\\[-1ex]\longrightarrow\end{array}
&
\left[\begin{array}{@{}cccc@{}}
\cdot&\cdot&\cdot&J\\
\cdot&\cdot&J&\cdot\\
\cdot&J&\cdot&\cdot\\
J&\cdot&\cdot&\cdot\\
\end{array}\right]
\end{array}$}
\vspace{-3.2ex}
\end{center}
\caption{The algorithm in the proof of Proposition~\ref{prop:chainEven} for $n=8$.
Here $J$ stands for $J_2$ and each dot represents a $2\times 2$ block of zeros.
The numbers above arrows indicate the length of a chain connecting the two matrices.
\label{fig:proofEven}}
\end{figure}

\begin{prop}
If $n\ge5$ is odd, then there is a chain of length $2n(n-2)-1$ from $P_n$ to $Q_n$ in the Bruhat order of $\mathcal{A}(n,2)$.
\label{prop:chainOdd}
\end{prop}

\begin{proof}
All chains in this proof are with respect to the Bruhat order of $(0,1)$--matrices.
Similarly to the proof of Proposition~\ref{prop:chainEven}, we give a 
construction of the desired chain.
For $n=5$ the chain of Lemma~\ref{lem:chain29} has the desired length.
Let $n=2k+5$ where $k\ge1$. 
Then $P_n=P_{2k}\oplus P_5$.
By applying the chain of Lemma~\ref{lem:chain29} to the submatrix of $P_n$
induced by its last five rows and its last five columns, we obtain a chain of
length $29$ from $P_n$ to $A=P_{2k}\oplus Q_5$.
Let $E$ be the $5\times5$ submatrix of $A$ induced by the rows $2k-1,2k,n-2,n-1,n$
and columns $2k-1,\ldots,2k+3$. Then $E=J_2\oplus F'_3$ and we may apply the
chain of Lemma~\ref{lem:chain24} to extend the current chain by $24$.
We repeat this procedure $k$ times to obtain a chain ending at 
\[
\def\arraystretch{1.4}
C=\left[\begin{array}{cc}0&P_{n-3}\\F'_3&0\end{array}\right].
\]
Using Proposition~\ref{prop:chainEven}, we may extend this chain to end at
\[
\def\arraystretch{1.4}
Q_n=\left[\begin{array}{cc}0&Q_{n-3}\\ F'_3&0\end{array}\right].
\]
This chain has length
\[29+24k+2(n-3)(n-5)=29+12(n-5)+2(n-3)(n-5)=2n(n-2)-1.\]
The algorithm presented in this proof is illustrated in Figure~\ref{fig:proofOdd} for $n=9$.
\end{proof}

\begin{figure}[ht]
\begin{center}
\scalebox{.96}{$\begin{array}{@{}c@{}c@{}c@{}c@{}c@{}c@{}c@{}c@{}c@{}}
\left[\begin{array}{@{}cccc@{}}
J&\cdot&\cdot&\cdot\\
\cdot&J&\cdot&\cdot\\
\cdot&\cdot&J&\cdot\\
\cdot&\cdot&\cdot&F\\
\end{array}\right]
&
\begin{array}{c}\scalebox{.8}{29}\\[-1ex]\longrightarrow\end{array}
&
\left[\begin{array}{@{}cccc@{}}
J&\cdot&\cdot&\cdot\\
\cdot&J&\cdot&\cdot\\
\cdot&\cdot&\cdot&J\\
\cdot&\cdot&F'&\cdot\\
\end{array}\right]
&
\begin{array}{c}\scalebox{.8}{24}\\[-1ex]\longrightarrow\end{array}
&
\left[\begin{array}{@{}cccc@{}}
J&\cdot&\cdot&\cdot\\
\cdot&\cdot&J&\cdot\\
\cdot&\cdot&\cdot&J\\
\cdot&F'&\cdot&\cdot\\
\end{array}\right]
&
\begin{array}{c}\scalebox{.8}{24}\\[-1ex]\longrightarrow\end{array}
&
\left[\begin{array}{@{}cccc@{}}
\cdot&J&\cdot&\cdot\\
\cdot&\cdot&J&\cdot\\
\cdot&\cdot&\cdot&J\\
F'&\cdot&\cdot&\cdot\\
\end{array}\right]
\\ \\[-1em]
\multicolumn{5}{r}{\small\text{(Using Proposition~\ref{prop:chainEven} for }n=6\text{)}}
& 
\begin{array}{c}\scalebox{.8}{48}\\[-1ex]\longrightarrow\end{array}
&
\left[\begin{array}{@{}cccc@{}}
\cdot&\cdot&\cdot&J\\
\cdot&\cdot&J&\cdot\\
\cdot&J&\cdot&\cdot\\
F'&\cdot&\cdot&\cdot\\
\end{array}\right]
\end{array}$}
\vspace{-3.2ex}
\end{center}
\caption{The algorithm in the proof of Proposition~\ref{prop:chainOdd} for $n=9$.
Here $J$ and $F$ stand for $J_2$ and $F_3$ respectively,
and each dot represents a block of zeros of the appropriate size.
The numbers above arrows indicate the length of a chain connecting the two matrices.\label{fig:proofOdd}}
\end{figure}

\begin{thm}
Let $n\ge4$ and let $\delta(n)$ denote the largest possible length of a chain in the Bruhat order of the class $\mathcal{A}(n,2)$.
Then
\[
\delta(n)=\begin{cases}
2n(n-2) & {\ \rm if\ }2|n,\\
2n(n-2)-1 & {\ \rm if\ }2\not|n.
\end{cases}
\]
\label{thm:maxlen}
\end{thm}

\begin{proof}
Let $A_0\ble A_1\ble \cdots \ble A_k$ be a chain in the Bruhat order of $\mathcal{A}(n,2)$.
By Lemma~\ref{lem:nu} we have $\nu(A_0)<\nu(A_1)<\cdots<\nu(A_k)$,
from which we obtain $k\le\nu(A_k)-\nu(A_0)$. Since a chain of maximum length $\delta(n)$ begins with a minimal element and ends with a maximal element, we obtain \[\delta(n)\le\max\set{\nu(Q)}-\min\set{\nu(P)},\]
where the maximum is over all maximal matrices $Q$ and
the minimum is over all minimal matrices $P$ in the Bruhat order of $\mathcal{A}(n,2)$.
On the other hand, from the characterization of minimal and maximal elements in the Bruhat order of $\mathcal{A}(n,2)$ discussed earlier in this section, since $\nu(J_2)=1$ and $\nu(F_3)=2$, a minimal matrix with smallest possible number of inversions cannot have more than one direct sum component of $F_3$. Therefore, $\nu(P_n)$ is the smallest possible value of $\nu(P)$ where $P\in\mathcal{A}(n,2)$ is minimal in the Bruhat order.
It can be shown similarly that $\nu(Q_n)$ is the largest possible value of $\nu(Q)$ where $Q\in\mathcal{A}(n,2)$ is maximal in the Bruhat order.
Therefore, we obtain
\[\delta(n)\le\nu(Q_n)-\nu(P_n).\]
Since there are no inversions in $P_n$ involving entries from different $J_2$ and $F_3$ direct sum components, we have $\nu(P_n)=\ceil{n/2}$.
In $Q_n$ on the other hand, every pair of ones in different $J_2$ and $F'_3$ blocks gives an inversion, while there are
$\nu(J_2)=1$ and $\nu(F'_3)=7$ inversions within each block.
A simple calculation gives $\nu(Q_n)=\floor{(4n^2-7n)/2}$.
We conclude that
\[\delta(n)\le\floor{(4n^2-7n)/2}-\ceil{n/2}=\begin{cases}
2n(n-2) & {\ \rm if\ }2|n,\\
2n(n-2)-1 & {\ \rm if\ }2\not|n.
\end{cases}\]
The constructions of Propositions~\ref{prop:chainEven} and~\ref{prop:chainOdd}
prove the lower bound.
\end{proof}

\begin{rem*}
Note that $\mathcal{A}({2,2})=\{J_2\}$ and $\mathcal{A}({3,2})$ is equivalent to $\mathcal{A}({3,1})$ via matrix complement. Therefore, $\delta(2)=0$ and $\delta(3)=3$. Of these,
only $\delta(2)$ agrees with the formula proved for $n\ge4$ in
Theorem~\ref{thm:maxlen}.
\end{rem*}

The chains found in Propositions~\ref{prop:chainEven} and~\ref{prop:chainOdd},
posses the property that the number of inversions in every two consecutive matrices
in a chain differ by~$1$. If $C,D\in\mathcal{A}(R,S)$, a chain 
\[C=A_0\sbleq A_1\sbleq A_2\sbleq\cdots\sbleq A_k=D\]
in the secondary Bruhat order
from $A$ to $C$ is called a {\em tight chain} if $\nu(A_i)=\nu(A_{i-1})+1$
for all $1\le i\le k$.
It is clear that such a tight chain (if one exists) has the maximum length among
all chains from $A$ to $C$ in the secondary Bruhat order.
A natural question to ask here is whether a tight chain exists between any two
matrices in $\mathcal{A}(n,2)$, $\mathcal{A}({n,k})$, or $\mathcal{A}(R,S)$ in general.

The techniques of this section can be applied to find the spectrum of lengths of
maximal chains in the Bruhat order of $\mathcal{A}(n,2)$. For $n\le5$, this spectrum consists
of only one length. On the other hand, for $n=6$, there are two minimal matrices
with $3$ and $4$ inversions, and two maximal matrices with $50$ and $51$ inversions.
Therefore, maximal chains in the Bruhat order of $\mathcal{A}({6,2})$ can have lengths $46$, $47$, or~$48$. Indeed tight chains of these lengths do exist between any pair of minimal and maximal matrices. This raises the following question.

\begin{pr}
Let $n\ge2$ and $A,C\in\mathcal{A}(n,2)$ such that $A\bleq C$.
Does there exist a tight chain in the (secondary) Bruhat order of $\mathcal{A}(n,2)$
from $A$ to $C$?
\label{q:tightchain}
\end{pr}

\begin{figure}[ht]
\begin{center}
\begin{tikzpicture}[scale=1.64]
\node [circle,draw, inner sep=1pt] (A1) at (-2.4142,0) {\scriptsize $A_1$};
\node [circle,draw, inner sep=1pt] (A2) at (-1,0) {\scriptsize $A_2$};
\node [circle,draw, inner sep=1pt] (A3) at (0,1) {\scriptsize $A_3$};
\node [circle,draw, inner sep=1pt] (A4) at (0,-1) {\scriptsize $A_4$};
\node [circle,draw, inner sep=1pt] (A5) at (1,0) {\scriptsize $A_5$};
\draw[->] (A1) -- (A2);
\draw[->] (A1) edge [bend left] (A3);
\draw[->] (A1) edge [bend right] (A4);
\draw[->] (A1) edge [bend right] (A5);
\draw[->] (A2) -- (A3);
\draw[->] (A2) -- (A4);
\draw[->] (A2) -- (A5);
\draw[->] (A3) -- (A5);
\draw[->] (A4) -- (A5);
\end{tikzpicture}
\vspace{-3.2ex}
\end{center}
\caption{The Bruhat graph of the class $\mathcal{A}({R,R})$ where $R=(2,2,1)$.
There is an arc from $A_i$ to $A_j$ when $A_i\bleq A_j$.\label{fig:A221}}
\end{figure}
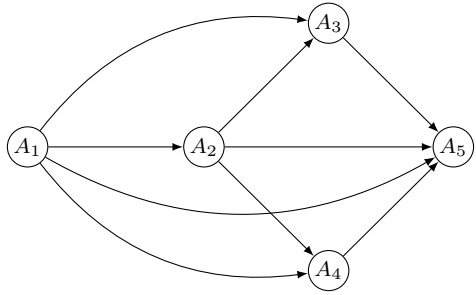

The statement in the above question does not hold for general classes $\mathcal{A}(R,S)$.
For example, if $R=S=(2,2,1)$,
then $\mathcal{A}(R,S)$ consists of the matrices
\begin{center}
$A_1=\left[\begin{array}{@{\ }c@{\ \ }c@{\ \ }c@{\ }}
1&1&0\\1&1&0\\0&0&1\end{array}\right]\text{, \ }
A_2=\left[\begin{array}{@{\ }c@{\ \ }c@{\ \ }c@{\ }}
1&1&0\\1&0&1\\0&1&0\end{array}\right]\text{, \ }
A_3=\left[\begin{array}{@{\ }c@{\ \ }c@{\ \ }c@{\ }}
1&1&0\\0&1&1\\1&0&0\end{array}\right]\text{,}$\\[2ex]
$A_4=\left[\begin{array}{@{\ }c@{\ \ }c@{\ \ }c@{\ }}
1&0&1\\1&1&0\\0&1&0\end{array}\right]\text{, \ and \ }
A_5=\left[\begin{array}{@{\ }c@{\ \ }c@{\ \ }c@{\ }}
0&1&1\\1&1&0\\1&0&0\end{array}\right].$
\end{center}
The Bruhat order on this class is illustrated in Figure~\ref{fig:A221}.
It can be seen in this graph that $A_5$ covers $A_4$ in the Bruhat order, thus
the maximum length of a chain from $A_4$ to $A_5$ is $1$, while
$\nu(A_5)-\nu(A_4)=6-3=3$. Note that $A_1$ and $A_5$
are the unique minimal and maximal elements of this class in the Bruhat order.
It is seen in Figure~\ref{fig:A221} that the maximum length of a chain from
$A_1$ to $A_5$ is $3$ while $\nu(A_5)-\nu(A_1)=5$.
Therefore, the above question does not hold for general classes $\mathcal{A}(R,S)$ even when
restricted to minimal and maximal matrices.

\section{Concluding remarks}

We define an inversion in a $(0,1)$--matrix.
Our main lemma proving that the number $\nu(A)$ of inversions in
a matrix $A\in\mathcal{A}(R,S)$ is monotonic with respect to the secondary
Bruhat order, provides a powerful tool in proving upper bounds
on the length of chains in this order. This feature is showcased
in our main result which establishes the maximum length of a chain
in the Bruhat order of the classes~$\mathcal{A}(n,2)$. This paper offers
several directions for future work, of which we have highlighted
two in Questions~\ref{q:ARSmonotone} and~\ref{q:tightchain}.

\section{Acknowledgements}

The author is grateful to Carlos da Fonseca who introduced him to
this fascinating problem.


\begin{thebibliography}{1}

\bibitem{Brualdi1980}
R.\,A. {Brualdi}.
\newblock {Matrices of zeros and ones with fixed row and column sum vectors}.
\newblock {\em {Linear Algebra Appl.}}, 33:159--231, 1980.

\bibitem{Brualdi2006}
R.\,A. {Brualdi}.
\newblock {Algorithms for constructing $(0,1)$-matrices with prescribed row and
  column sum vectors}.
\newblock {\em {Discrete Math.}}, 306(23):3054--3062, 2006.

\bibitem{BrualdiBook2006}
R.\,A. {Brualdi}.
\newblock {\em {Combinatorial matrix classes}}.
\newblock Number 108 in Encyclopedia of Mathematics and its Applications.
  Cambridge University Press, 2006.

\bibitem{BrualdiDeaett2007}
R.\,A. {Brualdi} and L. {Deaett}.
\newblock {More on the Bruhat order for $(0,1)$--matrices}.
\newblock {\em {Linear Algebra Appl.}}, 421(2-3):219--232, 2007.

\bibitem{BrualdiHwang2004}
R.\,A. {Brualdi} and S.-G. {Hwang}.
\newblock {A Bruhat order for the class of $(0,1)$-matrices with row sum vector
  $R$ and column sum vector $S$}.
\newblock {\em {Electron. J. Linear Algebra}}, 12:6--16, 2005.

\bibitem{ConflittiDaFonsecaMamede2012}
A. {Conflitti}, C.\,M. {da Fonseca}, and R. {Mamede}.
\newblock {The maximal length of a chain in the Bruhat order for a class of
  binary matrices}.
\newblock {\em {Linear Algebra Appl.}}, 436(3):753--757, 2012.

\bibitem{Ryser1957}
H.\,J. Ryser.
\newblock {Combinatorial properties of matrices of zeros and ones}.
\newblock {\em {Can. J. Math.}}, 9:371--377, 1957.

\bibitem{Ryser1963}
H.\,J. Ryser.
\newblock {\em {Combinatorial mathematics}}.
\newblock Number~14 in The Carus Mathematical Monographs. The Mathematical
  Association of America, 1963.

\end{thebibliography}

\def\soft#1{\leavevmode\setbox0=\hbox{h}\dimen7=\ht0\advance \dimen7
  by-1ex\relax\if t#1\relax\rlap{\raise.6\dimen7
  \hbox{\kern.3ex\char'47}}#1\relax\else\if T#1\relax
  \rlap{\raise.5\dimen7\hbox{\kern1.3ex\char'47}}#1\relax \else\if
  d#1\relax\rlap{\raise.5\dimen7\hbox{\kern.9ex \char'47}}#1\relax\else\if
  D#1\relax\rlap{\raise.5\dimen7 \hbox{\kern1.4ex\char'47}}#1\relax\else\if
  l#1\relax \rlap{\raise.5\dimen7\hbox{\kern.4ex\char'47}}#1\relax \else\if
  L#1\relax\rlap{\raise.5\dimen7\hbox{\kern.7ex
  \char'47}}#1\relax\else\message{accent \string\soft \space #1 not
  defined!}#1\relax\fi\fi\fi\fi\fi\fi}

\end{document}